\documentclass[leqno,a4paper,11pt]{article}
\usepackage{amsmath,amsthm,amssymb}
\input diagrams
\newcommand{\ar}{\ensuremath{\rightarrow}}
\newcommand{\Ar}{\ensuremath{\Rightarrow}}
\newcommand{\iar}{\ensuremath{\stackrel{\sim}{\longrightarrow}}}
\newcommand{\iso}{\ensuremath{\cong}}
\newcommand{\equi}{\ensuremath{\simeq}}
\newcommand{\Sets}{\ensuremath{\mathbf{Sets}}}
\newcommand{\E}{\ensuremath{\mathcal{E}}}
\newcommand{\mc}[1]{\ensuremath{\mathcal{#1}}}
\newcommand{\un}[1]{\lvert#1\rvert}
\newcommand{\iiff}{\quad\text{iff}\quad}
\newcommand{\iimp}{\quad\text{implies}\quad}
\newcommand{\aand}{\ \text{and}\ }
\newcommand{\abs}[1]{\lvert#1\rvert}
\newcommand{\mono}{\rightarrowtail}

\newcommand{\bb}[1]{\ensuremath{\mathbb{#1}}}
\DeclareMathOperator{\Log}{\mathbf{Log}}
\DeclareMathOperator{\Mod}{\mathbf{Mod}}
\DeclareMathOperator{\sh}{Sh}
\DeclareMathOperator{\true}{true}
\DeclareMathOperator{\false}{false}
\newtheorem{theorem}{Theorem}
\newtheorem{proposition}[theorem]{Proposition}
\newtheorem{lemma}[theorem]{Lemma}
\theoremstyle{remark}
\newtheorem{remark}[theorem]{Remark}
\theoremstyle{definition}
\newtheorem{definition}[theorem]{Definition}
\theoremstyle{plain}
\begin{document}
\title{Topological Completeness for Higher-Order Logic}
\author{%
	S.\ Awodey\thanks{%
		Philosophy Department, Carnegie Mellon University,
		Pittsburgh PA 15213-3890, USA. \tt{awodey@cmu.edu}}
	\and
	C.\ Butz\thanks{%
		BRICS, Basic Research in Computer Science, Centre of the
		Danish National Research Foundation, Computer Science
		Department, Aarhus University, Ny Munkegade, Bldg.\ 540, 8000
		Aahrus C., Denmark.  \tt{butz@brics.dk}}
	}
\maketitle

\begin{abstract}
\noindent Using recent results in topos theory, two systems of
higher-order logic are shown to be complete with respect to sheaf
models over topological spaces---so-called ``topological semantics''.
The first is classical higher-order logic, with relational
quantification of finitely high type; the second system is a
predicative fragment thereof with quantification over functions
between types, but not over arbitrary relations.  The second theorem
applies to intuitionistic as well as classical logic.
\end{abstract}


\section*{Introduction}

Higher-order logic (also known as ``type theory'') is logic that includes
quantification over functions or relations.  Many basic mathematical objects
and theories can only be defined using this logic; the natural numbers and
topological spaces are familiar examples.  A more precise specification of
what we call classical higher-order logic is given in \S 1 below.

As is well-known, higher-order theories are generally incomplete with
respect to (standard) models in \Sets; that is, $\bb{T}\models\sigma$
does not imply $\bb{T}\vdash\sigma$ for $\bb{T}$ a theory in
higher-order logic and $\vdash$ the entailment relation of any
reasonable deductive calculus.  It is by now also well-known that
higher-order logic can be modeled in suitable generalized categories
of sets, namely (elementary) topoi, and that with regard to such
topos-valued semantics, standard higher-order deduction is complete
(see \cite{L&S} for details).

Our results in this paper are concerned with topos models of a very
special and natural kind, namely sheaves over topological spaces.  If
$X$ is a space, a model in the category $\sh(X)$ of all sheaves on $X$
shall be called a \emph{topological model}.  We will show that
higher-order logic is complete with respect to such models; for the
reader unfamiliar with sheaf theory, we wish to emphasize their
elementary topological character.  Under the equivalence
$\sh(X)\equi\mathbf{Etale}/X$ a sheaf on a space $X$ is essentially
the same thing as an \'etale space over $X$:\ a space $E$ equipped
with a local homeomorphism $p : E\to X$ (also called an \'etale
bundle).  The various fibers $p^{-1}x$ of $E$ (the stalks of the
sheaf) for the points $x\in X$ may be regarded as sets varying
continuously over $X$.  A morphism of \'etale spaces is just a
continuous map $f : E\to E'$ over $X$, i.e.\ with $p'f=p$ as in the
commutative triangle %
\begin{diagram}
    E & \rTo^{f} &&& E' \\
    &\rdTo<p & & \ldTo>{p'}&\\
    && \ X. &&\\
    \end{diagram}
Products, exponentials (``function spaces''), etc.\ of \'etale spaces
of course agree with those calculated as sheaves.  A topological model
of a single-sorted theory thus consists of an \'etale space $p\colon
E\to X$ over a base space $X$ together with suitable operations, which
are simply continuous maps over~$X$.

As the reader who \emph{is} familiar with sheaf theory will have
noted, our topological models are just what are usually called
``sheaves of {\ldots}\!s'', at least in the case of equational,
first-order theories.  Thus a topological model of the theory of
groups is a sheaf of groups, and so on.

Despite the ultimately simple character of topological models, we use
the more general language and methods of sheaf theory and topoi to
study them.  Our first theorem, proved in \S 3 below, asserts the
completeness of standard, classical higher-order deduction
$\vdash^{c}$ with respect to such topological semantics.

\newtheorem*{theorem A}{Theorem A}
\begin{theorem A}
Let $\bb{T}$ be a higher-order theory.  There exists a classical topological
model $M$ of $\bb{T}$ such that, for any higher-order sentence $\sigma$ in the
language of $\bb{T}$,
\[
	\bb{T}\vdash^{c}\sigma\qquad\text{if and only if}\qquad
	M\models^{c}\sigma.
	\]
Moreover, the model $M$ has the property that every continuous function
between the interpretations of type symbols is logically definable.
\end{theorem A}

What permits theorem A to be true is our notion of a \emph{classical}
model.  In an arbitrary topological sheaf topos $\sh(X)$ there are two
natural candidates for the interpretation of the type $2$ of formulas
(or ``propositions'', or ``truth values'') of a higher-order theory;
to wit, the sheaf $\Omega$ of open subsets of $X$ and the coproduct
$1+1$.  In the language of \'etale spaces, $1+1$ is the double
covering $X\times 2\to X$.  As detailed in \S 2 below, a classical
model uses the latter to interpret the type of formulas.  Function and
power types are then interpreted as exponentials of sheaves (sometimes
called ``internal homs'' or ``sheaf-valued homs'').  This standard
treatment of exponentials is what chiefly distinguishes topological
models from so-called Henkin models (see the appendix below for the
exact relation between the two).  Thus in particular, for any type $Z$
the power type $2^{Z}$ is interpreted as the sheaf of
\emph{complemented} subsheaves of the interpretation of $Z$.  By
further requiring of a classical model that the types be interpreted
by so-called decidable sheaves, we can model classical higher-order
logic in non-boolean topoi like $\sh(X)$, which is impossible when
interpreting the type $2$ by the subobject classifier $\Omega$.
Indeed under that interpretation the analogue of theorem A
fails---even permitting arbitrary Grothendieck topoi in place of
topological sheaf topoi---as can be seen using G\"{o}del
incompleteness.

The issue of how to interpret the type of formulas of course vanishes when
one considers the fragment of higher-order logic that results from omitting
that type.  This fragment---which we call $\lambda$-logic and describe in
\S 4 below---may be regarded as a marriage of elementary logic and the
$\lambda$-calculus.  In addition to the usual propositional and
quantificational language of elementary logic, it includes equations
between and quantification over functions, functions of functions, etc.
But since there is no type of formulas, there is no quantification over
``propositional functions'', i.e.\ over relations.

Many familiar mathematical constructions, theorems, and proofs can be
formalized in $\lambda$-logic.  A simple example is Cayley's theorem that
every group is isomorphic to a group of permutations of its elements.  The
axiom of choice, in the familiar form
\[
	\forall x\in X\exists y\in Y.\varphi(x, y)%
		\Ar \exists f\in Y^{X}\forall x\in X.\varphi(x, fx),
	\]
is also a statement of $\lambda$-logic.  An example of a (non-elementary)
$\lambda$-theory is synthetic differential geometry, applications to which
of the present work shall be discussed elsewhere.

Our theorem B states the completeness of $\lambda$-logic with respect to
topological models.  More generally than theorem A, theorem B holds for
standard, intuitionistic deductive entailment $\vdash$.

\newtheorem*{theorem B}{Theorem B}
\begin{theorem B}
Let $\bb{T}$ be a $\lambda$-theory.  There exists a topological model $M$
of $\bb{T}$
such that, for any $\lambda$-sentence $\sigma$ in the language of $\bb{T}$,
\[
	\bb{T}\vdash\sigma\qquad\text{if and only if}\qquad M\models\sigma.
	\]
Moreover, the model $M$ has the property that every continuous function
between the interpretations of type symbols is logically definable.
\end{theorem B}

Theorem B rests more squarely on one of the main supports of theorem
A, namely a recent covering theorem for topoi due to the second author
and I.\ Moerdijk.  This covering theorem is the real heart of our
completeness theorems; we sketch its application to our situation as
an appendix to this paper.  So as not to obscure the conceptual
simplicity of this application, our treatment of the standard details
of higher-order syntax and topos semantics is held quite brief.

Before getting down to business, we make two remarks on the statements of
the completeness theorems.  First, each has the form ``there exists a model
$M$ such that $\bb{T}\vdash\sigma$ just if $M\models\sigma$'', rather than the
more familiar (for set-valued semantics) ``$\bb{T}\vdash\sigma$ just if
for all models $M$, $M\models\sigma$ ''.  The stronger form given here is
made possible by considering models in topoi other than \Sets.  The
situation is analogous to that of the familiar Heyting-valued
completeness theorem for first-order intuitionistic logic \cite{F&S}, which
follows directly from our theorem B and indeed is the inspiration thereof.
Second, and more substantially, the additional ``Moreover\ldots'' clause of
each theorem states a further property of the respective logical system
that may be termed ``definitional completeness''.  It ensures that any
function which is ``present in every model'' is logically definable.  As in
the case of deductive completeness, this definitional completeness is
established in a strong form simply by exhibiting a \emph{single} model in
which \emph{every} function of suitable type is definable.  In light of the
topological nature of the models at issue here, logical definability then
coincides with continuity in that ``minimal'' model.  For further
discussion of this property (in the context of the $\lambda$-calculus) we
refer to \cite{Awod:lambda}.

\paragraph*{Acknowledgments.} We have both benefitted greatly from
conversations with Ieke Moerdijk on the spatial covering theorem and
its logical applications.  The Stefan Banach Mathematical Research
Center in Warsaw, and the organizers of the Rasiowa memorial
conference held there in December 1996, are thanked for supporting our
collaboration.

\section{Theories in classical higher-order logic}

The systems of classical higher-order logic that we consider are
essentially the same as those presented in \cite{B&J, L&S}, which
in turn are modern formulations of \cite{Church:40}.  We summarize
one particular formulation for the reader's convenience and for the
special purposes of \S 4.

\emph{Type symbols} are built up inductively from a given list of
\emph{basic type symbols} $X_{1}, \ldots, X_{n}$ and the \emph{type of formulas}
$2$ by the \emph{type-forming operations} $Y\times Z$ and $Z^{Y}$.

\emph{Terms} are built up inductively from \emph{variables}
and a given list of \emph{basic terms} $c_{1}, \ldots, c_{m}$.  Each
variable and basic term has a type.  The terms and their types are as
follows, writing $\tau:Z$ for ``$\tau$ is a term of type $Z$''.
\begin{itemize}
	\item If $\tau_{1}:Z_{1}$ and $\tau_{2}:Z_{2}$, then
	$\langle\tau_{1},\tau_{2}\rangle:Z_{1}\times Z_{2}$.
	
	\item If $\tau:Z_{1}\times Z_{2}$, then $\pi_{1}(\tau):Z_{1}$
	and $\pi_{2}(\tau):Z_{2}$.
	
	\item If $\tau:Z$ and $y$ is a variable of type $Y$, then
	$\lambda y.\tau:Z^{Y}$.
	 	
	\item If $\alpha:Z^{Y}$ and $\tau:Y$, then $\alpha(\tau):Z$.
	
	\item If $\tau,\tau':Z$, then $\tau = \tau' : 2$.
	
	\item If $\varphi,\psi:2$ and $y$ is a variable of type $Y$, then the
	following are terms of type $2$:
	 \[
	 	\top,\ \bot,\ \neg\varphi,\ \varphi\wedge\psi,\ \varphi\vee\psi,\
	 	\varphi\Ar\psi,\ \forall y.\varphi,\ \exists y.\varphi.
	 	\]
\end{itemize}

A \emph{basic language} (or \emph{signature}) consists of
basic type symbols $X_{1}, \ldots, X_{n}$ and basic constant symbols
$c_{1}, \ldots, c_{m}$.  A \emph{theory} consists of a basic language and a
list of sentences (closed formulas) $\sigma_{1}, \ldots, \sigma_{k}$
therein, called \emph{axioms}.  Given a theory $\bb{T}$, the
\emph{language} $\mc{L}(\bb{T})$ of
$\bb{T}$ is the set of terms in the basic language of $\bb{T}$.

The \emph{entailment relation} $\varphi\vdash\psi$ between
formulas is specified in the usual way by a deductive calculus.  To include
the possibility of ``empty'' types, it is convenient to give a family of
entailment relations $\varphi\vdash_{\mathbf{x}}\psi$ indexed by lists
$\mathbf{x} = (x_{1},\ldots,x_{i})$ of distinct variables including all
those occurring free in $\varphi$ and $\psi$.  These relations are
generated by the following conditions (``rules of inference''):

\begin{enumerate}
  \item \emph{Order}\nopagebreak
  \begin{enumerate}
  	\item $\varphi\vdash_{\mathbf{x}}\varphi$
	\item $\varphi\vdash_{\mathbf{x}}\psi
		\aand\psi\vdash_{\mathbf{x}}\vartheta
		\iimp\varphi\vdash_{\mathbf{x}}\vartheta$	
	\item $\varphi\vdash_{\mathbf{x},y}\psi
		\iimp\varphi[\tau/y]\vdash_{\mathbf{x}}\psi[\tau/y]$			
	\end{enumerate}	
  \item \emph{Equality}
  \begin{enumerate}
  	\item $\top\vdash_{\mathbf{x}}\tau=\tau$
  	\item $\tau = \tau' \vdash_{\mathbf{x}} \varphi[\tau/y]\Ar\varphi[\tau'/y]$
	\item $\vartheta\vdash_{\mathbf{x}}\varphi\Ar\psi
			\aand\vartheta\vdash_{\mathbf{x}}\psi\Ar\varphi
		\iimp \vartheta\vdash_{\mathbf{x}}\varphi=\psi$
	\item $\forall y. \alpha(y)=\beta(y)\vdash_{\mathbf{x}} \alpha=\beta$
	\end{enumerate}
  \item \emph{Products}
  \begin{enumerate}
  	\item $\top\vdash_{\mathbf{x}}%
  		\langle\pi_{1}\tau,\pi_{2}\tau\rangle=\tau$
	\item $\top\vdash_{\mathbf{x}}%
		\pi_{i}\langle\tau_{1},\tau_{2}\rangle=\tau_{i},\quad i=1,2$
	\end{enumerate}	
  \item \emph{Exponents}
  \begin{enumerate}
  	\item $\top\vdash_{\mathbf{x}}(\lambda x.\tau)(x)=\tau$
	\item $\top\vdash_{\mathbf{x}}\lambda x.\alpha(x)=\alpha$
	\end{enumerate}
\pagebreak[2]		
  \item \emph{Elementary logic}
  \begin{enumerate}
  	\item $\bot\vdash_{\mathbf{x}}\varphi$
	\item $\varphi\vdash_{\mathbf{x}}\top$
	\item $\varphi\vdash_{\mathbf{x}}\neg\psi
		\iiff \varphi\wedge\psi\vdash_{\mathbf{x}}\bot$
	\item $\vartheta\vdash_{\mathbf{x}}\varphi
			\aand \vartheta\vdash_{\mathbf{x}}\psi
		\iiff \vartheta\vdash_{\mathbf{x}}\varphi\wedge\psi$
	\item $\vartheta\vee\varphi\vdash_{\mathbf{x}}\psi
		\iiff \vartheta\vdash_{\mathbf{x}}\psi
			\aand\varphi\vdash_{\mathbf{x}}\psi$
	\item $\vartheta\wedge\varphi\vdash_{\mathbf{x}}\psi
		\iiff \vartheta\vdash_{\mathbf{x}}\varphi\Ar\psi$
	\item $\vartheta\vdash_{\mathbf{x},y}\varphi
		\iiff \vartheta\vdash_{\mathbf{x}}\forall y.\varphi$
	\item $\exists y.\vartheta\vdash_{\mathbf{x}}\varphi
		\iiff \vartheta\vdash_{\mathbf{x},y}\varphi$
	\end{enumerate}	
  \end{enumerate}
In the foregoing, the $\tau$'s are arbitrary terms; $\varphi$,
$\psi$, $\vartheta$ are formulas; and $\alpha$, $\beta$ are terms
of the same exponential type.  In writing e.g.\
$\varphi[\tau/y]\vdash_{\mathbf{x}}\psi[\tau/y]$ in 1(c) it is
assumed that $\varphi[\tau/y]$ and $\psi[\tau/y]$ are formulas
with no free variables apart from $x_{1},\ldots,x_{i}$; so the
term $\tau$ must have the same type as the variable $y$ and no
other free variables.  As usual, the substitution notation
$\varphi[\tau/y]$ is understood to include a convention to avoid
binding free variables in $\tau$.

A sentence $\sigma$ is called \emph{provable} if $\top\vdash\sigma$, also
written $\vdash\sigma$.  For a theory \bb{T}, the notions of
\emph{$\bb{T}$-entailment} and \emph{$\bb{T}$-provability} are given by
adding the rules $\vdash\sigma$ for each axiom $\sigma$ of \bb{T}.

The \emph{classical entailment} relation $\vdash^{c}$ results from
$\vdash$ by adding the rule
\[
	\vdash^{c}\ \forall p.p\vee\neg p.
\]

\begin{remark}\label{succinct}
It is sometimes convenient to give a more succinct statement of the logical
calculus by defining some of the logical primitives in terms of others.  We
mention one particularly simple primitive basis which will be useful in the
next section (cf.\ \cite{L&S}).  Exponential types $Z^{Y}$ occur only in
the form $2^{Y}$ (``power types'', usually written $P(Y)$); $\lambda$-terms
$\lambda x.\varphi$ and evaluations $\alpha(\tau)$ are then restricted
accordingly, and more naturally written $\{x|\varphi\}$ and
$\tau\in\alpha$.  Projection operators $\pi_{i}(\tau)$ are eliminated in
favor of additional rules of inference.  The logical operations $\top,\
\bot,\ \neg,\ \wedge,\ \vee,\ {\Rightarrow},\ \forall,\ \exists$ are
defined in terms
of $=$ and $\langle -, - \rangle,\ \{x| - \},\ \in$.  We shall use the fact
that this primitive basis suffices in the following way:\ to interpret the
language of a theory it suffices to interpret the basic language, the type
of formulas, product and power types, and the term-forming operations
${\langle -, - \rangle},\ {\{x| - \}},\ {\in},\ {=}$.

In the opposite direction, one can enlarge the primitive logical basis by
including basic relation and function symbols in addition to basic
constant symbols, although these are not needed in the presence of higher
relation types.  Relation symbols will be useful in \S 4, however, where
there is no type of formulas; and both relation and function symbols are
used in elementary logic, where there are no higher types at all.
\end{remark}

\section{Semantics in topoi}

Let $\bb{T}$ be a theory in classical higher-order logic, as defined in the
foregoing section.  It is fairly obvious how to interpret $\bb{T}$ in an
arbitrary boolean topos \mc{B}: An \emph{interpretation} $M$ of $\bb{T}$ in
\mc{B} assigns to each basic type symbol $X$ an object $X_{M}$ of \mc{B},
and to
the type $2$ of formulas, the coproduct $1_{\mc{B}}+1_{\mc{B}}$ in \mc{B}
(which is
the subobject classifier),
\[
2_{M} = 1_{\mc{B}}+1_{\mc{B}}.
\]
The interpretation $M$ is then extended to product and power types by
setting
\begin{align*}
	(Y\times Z)_{M}& = Y_{M}\times Z_{M}&&\text{(product in \mc{B}),}\\
	(2^{Y})_{M}& = (2_{M})^{(Y_{M})}&&\text{(exponential in \mc{B}).}
	\end{align*}
On terms, $M$ assigns to each basic constant symbol $c$ of $\bb{T}$, having say
type $Z$, a morphism
\[
	c_{M} \colon  1_{\mc{B}}\ar Z_{M}
	\]
of \mc{B}, and variables are interpreted as identity morphisms.  The
interpretation is then extended inductively to all terms in
$\mc{L}(\bb{T})$ in the evident way, using the internal logic of
\mc{B} (cf.\ \cite[\S\S VI.5--7]{M&M}, also for the external meaning
of the logical operations thus modeled).  For example,
\[
	(\tau = \tau')_{M} = \delta\circ\langle\tau,\tau'\rangle_{M},
	\]
where $\delta\colon Z_{M}\times Z_{M}\ar 1_{\mc{B}}+1_{\mc{B}}$ classifies the
diagonal morphism $\Delta = \langle 1_{Z_{M}}, 1_{Z_{M}}\rangle \colon
Z_{M}\mono Z_{M}\times Z_{M}$, when $Z$ is the type of the terms $\tau,
\tau'$.

In particular, $M$ assigns to each formula $\varphi(y_{1}, \ldots, y_{n})$ with free
variables $y_{i}$ of types $Y_{i}$ a morphism
\[
	\varphi(y_{1}, \ldots, y_{n})_{M} \colon %
	 (Y_{1})_{M}\times\ldots\times (Y_{n})_{M}\longrightarrow 1_{\mc{B}}+1_{\mc{B}}
	 \]
of \mc{B}.  A sentence $\sigma$ is said to be \emph{true in $M$}, written
$M\models\sigma$, if
\[
 	\sigma_{M}=\true \colon  1_{\mc{B}}\ar 1_{\mc{B}}+1_{\mc{B}},
 	\]
where $\true \colon 1_{\mc{B}}\ar 1_{\mc{B}}+1_{\mc{B}}$ is the first
coproduct inclusion, which is the universal subobject.  Of course, an
interpretation $M$ is a \emph{model} of $\bb{T}$ if each axiom of $\bb{T}$
is true in $M$.

\subsection{Representing the category of models}\label{classtop}

\noindent Given models $M$ and $N$ of a theory $\bb{T}$ in a boolean topos
\mc{B}, there is an evident notion of an isomorphism $h\colon M\iar N$ of
$\bb{T}$-models, namely a family of isos $h=(h_{X}\colon X_{M}\iar X_{N})$
(indexed by the basic types $X$ of $\bb{T}$) that preserve the
interpretations of the constant symbols of $\bb{T}$, in the obvious sense.
Together with the evident composites and identities, one thus has for any
theory $\bb{T}$ and any boolean topos \mc{B} a \emph{category of models} of
$\bb{T}$ in \mc{B}, denoted
\[
	\Mod_{\bb{T}}(\mc{B}).
	\]
	
Observe that $\Mod_{\bb{T}}(\mc{B})$ is always a groupoid, i.e.\
a category in which every morphism is iso.  For example, if $\bb{T}$ is the
theory of topological spaces and \mc{B} is the topos \Sets, then
$\Mod_{\bb{T}}(\mc{B})$ is the category of all topological spaces and
homeomorphisms.  One can of course also consider other morphisms of models,
but the groupoid of isomorphisms suffices for our purposes.

A logical morphism between boolean topoi plainly preserves models and
their morphisms.  Such a functor $f \colon \mc{B}\ar\mc{B'}$ therefore
induces a functor
\[
	\Mod_{\bb{T}}(f)\colon \Mod_{\bb{T}}(\mc{B})\ar\Mod_{\bb{T}}(\mc{B'})
	\]
(a groupoid homomorphism) on the associated categories of models.

Now, every theory $\bb{T}$ in classical higher-order logic has a (higher-order)
\emph{classifying topos}, a boolean topos $\mc{B}_{\bb{T}}$
determined uniquely (up to equivalence) by the property: for any boolean
topos \mc{B} there is an equivalence of categories, natural in \mc{B},
\begin{equation}
	\Log(\mc{B}_{\bb{T}}, \mc{B})\ \equi\ \Mod_{\bb{T}}(\mc{B}),
	\label{univprop}
\end{equation}
where $\Log(\mc{B}_{\bb{T}}, \mc{B})$ is the category of logical morphisms
$\mc{B}_{\bb{T}}\ar\mc{B}$ and natural isomorphisms between them (cf.\
\cite{Awod:thesis}).  The classifying topos $\mc{B}_{\bb{T}}$ can be
constructed
``syntactically'' from $\mc{L}(\bb{T})$ in the style of \cite{B&J, L&S}; in
particular, it is a small category (indeed, it is countable).  In virtue of
its universal mapping property \eqref{univprop}, $\mc{B}_{\bb{T}}$ is freely
generated as a boolean topos by the ``universal model'' $U_{\bb{T}}\in
\Mod_{\bb{T}}(\mc{B}_{\bb{T}})$ associated to the identity logical morphism
$\mc{B}_{\bb{T}}\ar\mc{B}_{\bb{T}}$ under \eqref{univprop}.  By the syntactic
construction of $\mc{B}_{\bb{T}}$ this universal model has the following
properties, which we record for later use:

\begin{proposition}\label{univmodel}
\begin{enumerate}
	\item[(i)] For any sentence $\sigma\in\mc{L}(\bb{T})$,
	\[
		\bb{T}\vdash^{c}\sigma\quad\text{just if}\quad U_{\bb{T}}\models\sigma.
		\]
	\item[(ii)] For any types $Y$ and $Z$ and any morphism $f \colon
	Y_{U_{\bb{T}}}\ar Z_{U_{\bb{T}}}$ in $\mc{B}_{\bb{T}}$, there is a formula
	$\varphi(y,z)\in \mc{L}(\bb{T})$ such that
	\[
		\mathrm{graph}(f)=\{\langle y,z\rangle \mid \varphi(y,z)\}_{U_{\bb{T}}}
		\]
	(as subobjects of $Y_{U_{\bb{T}}}\times Z_{U_{\bb{T}}}$).
	\end{enumerate}
\end{proposition}

Observe that (i) of proposition \ref{univmodel} and the universal mapping
property \eqref{univprop} together entail the soundness and completeness of
the deductive calculus of \S 1 with respect to topos semantics:\
$\bb{T}\vdash^{c}\sigma$ if and only if for every $\bb{T}$-model
$M$, $M\models\sigma$.

\subsection{Classical semantics}

We now extend the foregoing topos semantics for classical higher-order
logic to non-boolean topoi.  Let $\bb{T}$ be a fixed theory and \mc{E} an
arbitrary topos.  We begin with a bit of notation: Let $\true\colon
1_{\mc{E}}\to\Omega_{\mc{E}}$ be the subobject classifier in \mc{E}, and
let us write
\[
     \abs{-} = (\true,\false) \colon
1_{\mc{E}}+1_{\mc{E}}\longrightarrow\Omega_{\mc{E}}
     \]
for the canonical map from the coproduct which, observe, is a monomorphism.
An arbitrary morphism $\varphi \colon  E\ar\Omega_{\mc{E}}$ of \mc{E} factors
through $\abs{-}$ just if the subobject $S_{\varphi}\mono E$ it classifies is
\emph{complemented}, i.e.\ if there is a subobject $S\mono E$ with
$S_{\varphi}+S\iso E$ (canonically).  When this is the case, let us write
$\overline{\varphi} \colon  E\ar 1_{\mc{E}}+1_{\mc{E}}$ for the unique morphism
such that
\[
	\varphi = \abs{\overline{\varphi}},
	\]
as indicated in
\begin{diagram}
     E &\rDashto<{\overline{\varphi}}& 1_{\mc{E}}+1_{\mc{E}}\\
     & \rdTo_{\varphi} & \dTo>{\ \abs{-}} \\
     && \Omega_{\mc{E}}.\\
    \end{diagram}
Recall that an object $E$ of \mc{E} is said to be \emph{decidable} if its
diagonal $\Delta \colon E\mono E\times E$ is complemented, thus just if
$\overline{\delta} \colon E\times E \ar 1_{\mc{E}}+1_{\mc{E}}$ exists.

Next, we define an \emph{interpretation} of the basic language of $\bb{T}$ in
\mc{E} exactly as in a boolean topos; in particular the type $2$ of
formulas is interpreted as $1_{\mc{E}}+1_{\mc{E}}$, which is plainly
decidable.  An interpretation $M$ such that for each type symbol $Z$ the
object $Z_{M}$ in \mc{E} is decidable shall be called a \emph{classical
interpretation} (or \emph{c-interpretation}).

Finally, by remark \ref{succinct} any c-interpretation $M$ can be extended
to all of
$\mc{L}(\bb{T})$ exactly as in a boolean topos, by interpreting the
term-forming
operations $\langle -, - \rangle,\ \{x| - \},\ \in$ as before and taking
$\overline{\delta}\colon Z_{M}\times Z_{M}\ar 1_{\mc{E}}+1_{\mc{E}}$ to
interpret $=$ at each type $Z$.  Thus just as before a c-interpretation $M$
assigns to each formula $\varphi(y_{1}, \ldots, y_{n})$ with free variables
$y_{i}$ of types $Y_{i}$ a morphism
\[
	\varphi(y_{1}, \ldots, y_{n})_{M} \colon %
	 (Y_{1})_{M}\times\ldots\times (Y_{n})_{M}%
	 	\longrightarrow 1_{\mc{E}}+1_{\mc{E}}.
	 \]
\begin{definition}
The relation $\models^{c}$ of \emph{satisfaction} for c-interpretations is
defined by:
\[
	M\models^{c}\sigma\iiff \abs{\sigma_{M}} = \true.
	\]
Thus a c-interpretation $M$ satisfies a sentence $\sigma$ just if the
following triangle commutes
\begin{diagram}
     1_{\mc{E}} &\rTo<{\sigma_{M}} & 1_{\mc{E}}+1_{\mc{E}}\\
     & \rdTo_{\true} & \dTo>{\ \abs{-}} \\
     && \Omega_{\mc{E}}.\\
    \end{diagram}
\end{definition}

A \emph{c-model} of the theory $\bb{T}$ is of course a
c-interpretation that satisfies the axioms of $\bb{T}$.  A
c-interpretation $M$ is therefore a c-model just if for each axiom
$\sigma$ the interpretation
\[
     \sigma_{M}\colon 1_{\mc{E}} \longrightarrow 1_{\mc{E}}+1_{\mc{E}}
     \]
is the first coproduct inclusion, just as in the boolean case.  Indeed, if
the topos \mc{E} is boolean, then every object is decidable, and a c-model
in \mc{E} is the same thing as a model.

\begin{proposition}[Soundness]
If $M$ is a c-model then for any sentence $\sigma$,
\[
	\bb{T}\vdash^{c}\sigma\iimp M\models^{c}\sigma.
	\]
\end{proposition}

\begin{proof}
Consider the classifying topos $\mc{B}_{\bb{T}}$ with universal model
$U_{\bb{T}}$.  There is an evident functor $m : \mc{B}_{\bb{T}}\to\E$
with $M=m(U_{\bb{T}})$ and
\[
	\sigma_{M}=\sigma_{m(U_{\bb{T}})}=m(\sigma_{U_{\bb{T}}})
	\]
for each sentence $\sigma$.  Although $m$ is not logical if \E\ is not
boolean, it still takes $\true : 1_{\mc{B}_{\bb{T}}}\ar
1_{\mc{B}_{\bb{T}}}+1_{\mc{B}_{\bb{T}}}$ to $\overline{\true} :
1_{\mc{E}}\ar 1_{\mc{E}}+1_{\mc{E}}$.  The claim thus follows from the
soundness of standard topos semantics (in particular, from proposition
\ref{univmodel}).
\end{proof}

\begin{remark}
If the interpretation $Z_{M}$ is decidable then for any type $Y$ the
canonical inclusion $(Z_{M})^{(Y_{M})}\mono \Omega^{(Y_{M}\times
Z_{M})}$ factors as indicated in the following diagram.
\begin{diagram}
     && (1_{\mc{E}}+1_{\mc{E}})^{(Y_{M}\times Z_{M})}\\
     & \ruTo & \dTo>{\ \abs{-}^{(Y_{M}\times Z_{M})}} \\
     (Z_{M})^{(Y_{M})} &\rEmbed & \Omega^{(Y_{M}\times Z_{M})}\\
    \end{diagram}
Thus even when defined in terms of power types as mentioned in remark
\ref{succinct}, the exponential types $Z^{Y}$ are still interpreted as
exponentials by a c-interpretation.
\end{remark}

\section{Topological completeness}

In this section we consider small topoi equipped with the finite epi
topology.  The covering families for this Grothendieck topology on a
topos \E\ are finite families of morphisms $(C_{i}\ar E)_{i}$ such
that the canonical map $\coprod_{i}C_{i}\ar E$ is epic.  Of the
following two technical lemmas, we omit the straightforward proof of
the first; its second part is folklore.

\begin{lemma}\label{lem:ypreserves}
\begin{enumerate}
  	\item[(i)] The finite epi topology is subcanonical.
	
	\item[(ii)] For each morphism $e : E'\to E$ in \E, the sheafified
	Yoneda embedding $y:\E\to\sh(\E)$ preserves not only the
	pullback functor $e^{*} : \E/E\to \E/E'$, but also its left
	and right adjoints,
	\[
   		\Sigma_{e}\dashv e^{*}\dashv\Pi_{e}  : \E/E'\to \E/E.
   		\]
	(Indeed, this is true for any subcanonical topology on a small
	category and any locally cartesian closed structure present
	there.)
  	\end{enumerate}
\end{lemma}

\begin{lemma}\label{lem:1}
Let $F\colon \mc{B} \ar\mc{E} $ be a left-exact functor from a
boolean topos \mc{B} to any topos \mc{E}.  If $F$ is continuous
for the finite epi topology then it preserves finite coproducts
and first-order logic.  If $F$ also preserves exponentials, then
it preserves c-models.
\end{lemma}

\begin{proof}
An object of a topos has an empty covering family for the finite epi
topology just if it is initial; so the continuous functor $F$
preserves initial objects.  The coproduct inclusions $B_{1},
B_{2}\mono B_{1}+B_{2}$ are a covering family of monos with
$B_{1}\wedge B_{2}=0\mono B_{1}+B_{2}$.  Since $F$ is also left-exact,
it then preserves coproducts as well.  Moreover, it then preserves boolean
complements of subobjects, whence it preserves negation $\neg$ since
\mc{B} is boolean.  The logical operations $\Rightarrow$ and $\forall$
are then also preserved, since in a boolean topos these can be
constructed from negation and operations that are preserved by
left-exact, continuous functors generally.  Finally, if
$F$ also preserves exponentials then by the foregoing it
preserves the interpretations of all types and the associated
term-forming operations, in addition to first-order
logic; whence it clearly preserves c-models.
\end{proof}

Theorem A will now follow by applying the covering theorem of the
appendix, which states that every Grothendieck topos with enough
points can be covered by a topological space via a connected, locally
connected geometric morphism.  We remind the reader that a
Grothendieck topos $\mc{G}$ is said to have \emph{enough points} if
the geometric morphisms $p : \Sets\to\mc{G}$ are jointly surjective
(i.e.\ if the inverse images $p^{*} : \mc{G}\to\Sets$ of these are
jointly faithful), and that a geometric morphism $f^{*}\dashv f_{*}
\colon \mc{E}\ar\mc{F}$ of topoi is \emph{connected} if the inverse
image functor $f^{*}$ is full and faithful, and \emph{locally
connected} (\cite{Barr&Pare}: ``molecular'') if $f^{*}$ commutes with
$\Pi$-functors.

\begin{theorem A}
Let $\bb{T}$ be a higher-order theory.  There exists a topological space
$X_{\bb{T}}$
and a c-model $M$ of $\bb{T}$ in $\sh(X_{\bb{T}})$ such that:
\begin{enumerate}
	\item[(i)] for any sentence $\sigma\in \mc{L}(\bb{T})$,
	\[
		\bb{T}\vdash^{c}\sigma\qquad\text{if and only if}\qquad M\models^{c}\sigma;
		\]
	\item[(ii)] given types $Y, Z$, every continuous function $f \colon  Y_{M}\ar
	Z_{M}$ over $X_{\bb{T}}$ is definable: there is a formula $\varphi(y,z)\in
	\mc{L}(\bb{T})$ such that
	\[
		\mathrm{graph}(f)=\{\langle y,z\rangle | \varphi(y,z)\}_{M}
		\]
	(as subsheaves of $Y_{M}\times Z_{M}$).
	\end{enumerate}
\end{theorem A}

\begin{proof}
First, one has the universal model $U_{\bb{T}}$ in the classifying
topos $\mc{B}_{\bb{T}}$, as in \S 2.1.  The Grothendieck topos
$\sh(\mc{B}_{\bb{T}})$ of sheaves on $\mc{B}_{\bb{T}}$ for the finite
epi topology is coherent, and so has enough points (cf.\ \cite{M&M}).
The covering theorem of the appendix therefore guarantees the
existence of a topological space $X_{\bb{T}}$ and a connected,
locally-connected geometric morphism
\[
	m \colon   \sh(X_{\bb{T}})\ar \sh(\mc{B}_{\bb{T}}).
	\]
The inverse image $m^{*}\colon \sh(\mc{B}_{\bb{T}})\ar \sh(X_{\bb{T}})$ of
$m$ satisfies all hypotheses of the foregoing lemma \ref{lem:1}, as does
the sheafified Yoneda embedding
\[
	y \colon   \mc{B}_{\bb{T}}\ar \sh(\mc{B}_{\bb{T}}).
	\]	
In particular, these functors preserve exponentials since they
preserve $\Pi$-functors (using lemma \ref{lem:ypreserves}).  The
composite $m^{*}\circ y\colon \mc{B}_{\bb{T}}\ar \sh(X_{\bb{T}})$
therefore also satisfies the hypotheses of lemma \ref{lem:1}, whence
one has the c-model
\[
	M = m^{*}\circ y(U_{\bb{T}})
	\]
in $\sh(X_{\bb{T}})$.  Since each of its factors is full and faithful, so
is the
functor $m^{*}\circ y$; the assertions (i) and (ii) thus follow from
proposition \ref{univmodel}.
\end{proof}

\begin{remark}(Infinitary generalizations) Theorem A clearly applies
equally to ``theories'' $\bb{T}$ with infinitely many type and/or constant
symbols and/or axioms, since in such cases the foregoing proof can begin
with a small topos $\mc{B}_{\bb{T}}$ which is a suitable colimit of
classifying topoi for (finite) theories.  We also merely mention that for
the case of infinitary logic, with set-indexed meets and joins of formulas,
a theorem analogous to theorem A holds, with complete Heyting algebras in
place of topological spaces.
\end{remark}

\section{$\lambda$-logic}

What we call $\lambda$-logic differs from classical higher-order logic
in that it has no type $2$ of formulas.  Type symbols are now built up
inductively from basic type symbols by the operations $-\times ?$ and
$-^{?}$.  Terms are built up inductively from variables, basic
constant symbols, and just the term-forming operations $\langle
-,?\rangle$, $\pi_{1}(-)$, $\pi_{2}(-)$, $\lambda y.(-)$, and $?(-)$.
Formulas are then constructed from terms and basic relation symbols in
the customary way, using the language of first-order logic with
equality.  Finally, a \emph{$\lambda$-theory} consists of (finitely
many) basic type, constant, and relation symbols, and closed formulas
in these parameters.

As rules of inference for the (intuitionistic) entailment relation
$\varphi\vdash_{\mathbf{x}}\psi$ on formulas one may take a standard deductive
calculus for (intuitionistic) many-sorted, first-order logic with
equality, augmented by the usual rules for the (typed)
$\lambda$-calculus.  Indeed, the rules given in \S 1 above are
suitable, under the omission of 2(c).

The notion of a model of a $\lambda$-theory in a topos is essentially
the same as that already given in \S 2.  It is, however, now more
natural to interpret basic relation symbols and other formulas by
subobjects (rather than their classifying morphisms), as is usually
done for first-order logic (cf.\ \cite{M&M}).  In particular, the
equality sign $=$ is interpreted in the standard way as a diagonal
morphism, and since classical logic is not being assumed, the notion
of a c-model is not required.

Deduction is clearly sound with respect to such semantics.  To show
that it is also complete---even with regard to just topological
models---one can proceed as in the classical higher-order case in \S3:
\begin{enumerate}
\item[(i)] Construct the syntactic category $\mc{S}_{\bb{T}}$ of
provable equivalence classes of formulas, to be equipped with the
finite epi topology (which is sub-canonical).
\item[(ii)] Apply the sheafified Yoneda embedding $y \colon
\mc{S}_{\bb{T}}\ar\sh(\mc{S}_{\bb{T}})$ (which preserves
$\lambda$-logic by lemma \ref{lem:ypreserves}) to get a full and
faithful model in a Grothen\-dieck topos with enough points.
\item[(iii)] Apply the covering theorem of the appendix to get a
connected, locally connected geometric covering map
$\sh(X_{\bb{T}})\ar\sh(\mc{S}_{\bb{T}})$ from a topological sheaf
topos $\sh(X_{\bb{T}})$.
\end{enumerate}
We leave it to the reader to fill in the details of this sketch to
provide the proof of the following.

\begin{theorem B}
Let $\bb{T}$ be a $\lambda$-theory.  There exists a topological space
$X_{\bb{T}}$
and a model $M$ of $\bb{T}$ in $\sh(X_{\bb{T}})$ such that:
\begin{enumerate}
	\item[(i)] for any $\lambda$-sentence $\sigma$ in the language of $\bb{T}$,
	\[
		\bb{T}\vdash\sigma\qquad\text{if and only if}\qquad M\models\sigma;
		\]

	\item[(ii)] given types $Y, Z$, every continuous function $f \colon  Y_{M}\ar
	Z_{M}$ over $X_{\bb{T}}$ is definable: there is a $\lambda$-formula
	$\varphi(y,z)$ in the language of $\bb{T}$ such that
	\[
		\mathrm{graph}(f)=\{\langle y,z\rangle | \varphi(y,z)\}_{M}
		\]
	(as subsheaves of $Y_{M}\times Z_{M}$).
	\end{enumerate}
\end{theorem B}

\section*{Appendix: The spatial cover}

In the proofs of theorems~A and~B, use was made of the following covering
theorem for topoi, which is part of theorem~13.5 of~\cite{Butz:thesis}
(also see \cite{Butz&Moer}; cf.\ \cite{J&M} for a related result).

\newtheorem*{cover}{Covering theorem}
\begin{cover}
For any Grothendieck topos $\cal{G}$ with enough points there is a
topological space $X_{\cal{G}}$ and a connected, locally connected
geometric morphism
\[
	\phi\colon\sh(X_{\cal{G}})\to {\cal G}.
	\]
\end{cover}

\noindent Thus in particular the inverse image functor $\phi^*\colon{\cal
G}\to\sh(X_{\cal G})$ is fully faithful and preserves exponentials and the
internal first-order logic of ${\cal G}$.

The purpose of this appendix is to describe the space $X_{\cal{G}}$ and the
covering map $ \phi \colon\sh(X_{\cal G})\to {\cal G}$ in the case of
principal interest here, namely when ${\cal G}=\sh({\cal B}_\bb{T})$ for ${\cal
B}_\bb{T}$ the small classifying topos of a (classical) higher-order theory,
equipped with the finite epi topology.  Thus we consider the situation of
theorem~A; that of theorem~B of course has a similar description.  Before
going into details, let us mention that in fact there are many different
spaces which will do the job, depending on various parameters that one is
free to choose.  We exhibit here just one such choice, intended to be
illuminating.

To begin, recall from \cite{Hen50} that classical higher-order logic
is complete with respect to {\em general models\/}, nowadays called
{\em Henkin models\/}.  The basic feature of a Henkin model $M$ of a
theory $\bb{T}$ is that a function type $Z^Y$ (or power type~$2^Y$) is
interpreted by a {\em subset\/} $(Z^Y)_M\subset (Z_M)^{(Y_{M})}$ of
the set of all functions from $Y_M$ to $Z_M$ (resp.\ of the power set
$\wp Y_{M}$), rather than by the set itself.  Of course, certain
closure conditions also have to be satisfied.  We mention only by the
way that such models can be shown to arise ``naturally'' as images of
the universal model $U_{\bb{T}}$ under continuous, left exact functors
$\mathcal{B}_{\bb{T}}\to\mathbf{Sets}$, and that the said completeness
can be inferred from this fact.  For the following, it will be convenient
to define the underlying set or \emph{universe} $\un{M}$ of a Henkin
model $M$ to be the (disjoint) union of the sets $Z_{M}$ for all types $Z$,
\[
	\un{M} =\bigcup\{Z_M\mid\mbox{$Z$ a type}\}.
	\]
	
To define the space $X_\bb{T}$ for the topos $\sh({\cal B}_\bb{T})$, fix a
\emph{sufficient} set $S_\bb{T}$ of countable Henkin models $M$ of
$\bb{T}$, i.e.\
$S_\bb{T}$ satisfies:
\[
	\text{$M\models\sigma$ for all $M\in S_\bb{T}$}
		\qquad\text{implies}\qquad\bb{T}\vdash\sigma
\]
for all $\bb{T}$-sentences $\sigma$.  For example, we could take (a set of
representatives of) all countable Henkin models of $\bb{T}$ as the set~$S_\bb{T}$.
We then define a \emph{labeling} of a Henkin model $M$ in $S_\bb{T}$ to be a
partial function
\[
	\mathbb{N}\supset \mathrm{dom}(\alpha)%
		\stackrel{\alpha}{\longrightarrow}\un{M}
\]
such that for each $a\in\un{M}$ the fiber $\alpha^{-1}(a)$ is infinite.

The points of the space $X_\bb{T}$ are labeled Henkin models in $S_\bb{T}$,
i.e.\
pairs
\[
	(M,\alpha)
\]
where $M\in S_\bb{T}$ and $\alpha$ is a labeling of~$M$.  The topology is
generated by basic open sets of the form
\[
	U_{\varphi(\bar{z}),\bar{n}}=%
		\{(M,\alpha)\mid
			\begin{array}[t]{l}
			\mbox{$\alpha(n_i)$ is defined and of type $Z_i$,}\\
			\mbox{and $M\models\varphi(\alpha(n_1),\ldots,\alpha(n_m))$}
				\quad\}
			\end{array}
\]
for $\varphi(\bar{z})=\varphi(z_1,\ldots,z_m)$ a $\bb{T}$-formula and
$\bar{n}=(n_1,\ldots,n_m)$ a tuple of natural numbers.

To describe the covering map $ \phi \colon\sh(X_{\bb{T}})\to \sh({\cal
B}_\bb{T})$ we
sketch the construction of the c-model $\Phi$ in $\sh(X_{\bb{T}})$ induced
by~$\phi^*$.  Here we use the equivalence, mentioned in the introduction,
$\sh(X_\bb{T})\simeq{\rm Etale}/X_\bb{T}$ of sheaves on $X_\bb{T}$ and \'etale
bundles over~$X_\bb{T}$.  For each type $Z$ we have the set
$$Z_\Phi=\sum_{(M,\alpha)\in X_\bb{T}}Z_M,$$
with the evident projection
$$\pi_Z\colon Z_\Phi\to X_\bb{T}.$$
We generate a topology on $Z_\Phi$ by declaring to be open:
\begin{itemize}
\item the sets $\pi^{-1}_Z(U)$ for $U\subset X_\bb{T}$ open (thus making
$\pi_Z$
continuous),
\item the sets $V_n=\{(M,\alpha,a)\mid\mbox{$a\in Z_M$,
$\alpha(n)$ is defined, and
$\alpha(n)=a$}\}$.
\end{itemize}

It is easily checked that $\pi_Z$ then becomes a local homeomorphism
(an \'etale map).  The assignment $Z\mapsto Z_{\Phi}$ extends in the
obvious way to a left exact, continuous functor ${\cal
B}_\bb{T}\to\sh(X_{\bb{T}})$ that preserves exponentials, inducing the
covering map $ \phi \colon\sh(X_{\bb{T}})\to \sh({\cal B}_\bb{T})$.
Finally, the stalk $x^*\Phi$ of the c-model $\Phi$ at a point
$x=(M,\alpha)$ of $X_\bb{T}$ is just the Henkin model $M$ itself,
which gives the relationship between our results and \cite{Hen50}.

\providecommand{\bysame}{\leavevmode\hbox to3em{\hrulefill}\thinspace}

\end{document}